\documentclass[12pt]{amsart}

\numberwithin{equation}{section}

\def\Im{{\text{\rm Im}}}
\def\CR{{\text{\rm CR}}}
\def\contr{\text{\rm contr}}
\def\pa{\partial}

\def\wt{\widetilde}
\def\wh{\widehat}
\def\MA{Monge-Amp\`ere }
\def\E{{\mathcal{E}}}

\def\C{\mathbb C}
\def\R{\mathbb R}

\def\CR{{\operatorname{CR}}}
\def\G{\mathcal{G}}

\def\Up{\Upsilon}
\def\up{\Up}

\newtheorem{theorem}{Theorem}[section]
\newtheorem{lemma}[theorem]{Lemma}
\newtheorem{proposition}[theorem]{Proposition}

\numberwithin{equation}{section}

\begin{document}

\title[Ambient metric construction of $Q$-curvature]%
{Ambient metric construction of $Q$-curvature
in conformal and CR geometries}

\author[C. Fefferman]{Charles Fefferman}
\address{
Department of Mathematics\\
Princeton University\\
Princeton, NJ 08544, USA
}
\email{cf@math.princeton.edu}

\author[K. Hirachi]{Kengo Hirachi
}
\address{
Graduate School of Mathematical Sciences\\
The University of Tokyo\\
Komaba, Megro, Tokyo 153-8914, Japan
}
\email{hirachi@ms.u-tokyo.ac.jp}


\thanks{The work of the second author was supported by Grant-in-Aid  for Scientific Research, JSPS}

\maketitle

\section{Introduction}
This article presents a geometric derivation of the
$Q$-curvature in terms of
the ambient metric associated with conformal and CR structures.
The $Q$-curvature in conformal geometry is a scalar Riemannian
invariant $Q$ that is conformally invariant up to an error
given by a conformally invariant power of the Laplacian.
In dimension $2$, the $Q$-curvature is
the half of the scalar curvature, $Q_1=R/2$, which satisfies
$$
e^{2\up}\wh Q_1=Q_1+\Delta\up\quad\text{whenever }\wh g=e^{2\up}g.
$$
In general even dimension $n$, the Laplacian $\Delta$ is replaced
by a conformally invariant $(n/2)^\text{th}$ power of the Laplacian
$P_{n/2}$ and the $Q$-curvature, $Q_{n/2}$, is
$\frac{1}{2(n-1)}\Delta^{n/2}R$ modulo nonlinear terms in curvature.
We here give a simple formula for $Q_{n/2}$ which directly follows
from the ambient metric construction of $P_{n/2}$ given in \cite{GJMS};
this formula can be generalized to any invariant differential operators $P$
on functions (densities of weight $0$) that arise in the ambient metric
construction. We also apply the construction of the $Q$-curvature to CR
geometry; it then turns out that the $Q$-curvature gives the coefficient
of the logarithmic singularity of the Szeg\"o kernel of 3-dimensional CR
manifolds.

The $Q$-curvature in general even dimensions was first defined by Branson
\cite{Br} in a study of the functional determinant of the conformal
Laplacian.  He used an argument of analytic continuation in the dimension,
in which the $Q$-curvature in dimension $n$ is defined from $P_{n/2}$ in
dimension $m>n$. For example, $Q_1=R/2$ in dimension $2$ is obtained from
the zeroth order term of conformal Laplacian $P_1=\Delta+\frac{(n-2)}{4(n-1)}R$
in dimension $n>2$. It was thus natural to ask: what is $Q_{n/2}$ in $n$
dimensional conformal geometry? An answer to this question was given by
 Graham-Zworski \cite{GZ} in their study of the scattering for the Laplacian
$\Delta_+$ in the Poincar\'e metric associated with conformal structure
$(M,[g])$ on the boundary at infinity.  They gave a formula for the
$Q$-curvature in terms of the scattering matrix, where the argument of
analytic continuation in dimension is replaced by the analytic continuation
in a spectral parameter. Their formulation of the $Q$-curvature was
significantly simplified in Fefferman-Graham \cite{FG2}; the $Q$-curvature
is given as the logarithmic term in the formal solution to a Dirichlet problem
for $\Delta_+$.  This construction has an intimate relation to the volume
expansion for the a Poincar\'e metric \cite{GZ}; in particular, it is shown
that the integral of $Q$ gives the coefficient of the log term of the expansion.

Our approach is directly related to the original derivation of the operator
$P_{n/2}$ of \cite{GJMS} in terms of the ambient metric space, a formally
constructed $(n+2)$-dimensional pseudo-Riemannian space $\wt\G$ that contains
the metric bundle $\G=\{t^2g\in S^2T^*M:t\in C^\infty(M)\}$ as a hypersurface.
The invariant differential operator $P_{n/2}$ arises in two ways:

\smallskip

(a) as an operator on the homogeneous functions on $\G$ induced from the
powers of the Laplacian $\wt\Delta^{n/2}$ on the ambient space;

\smallskip

(b) as an obstruction to the existence of a smooth homogeneous solution to
$\wt\Delta F=0$ with an initial condition on $\G$.

\smallskip

\noindent
Corresponding to these derivations, we give two formulas for the
$Q$-curvature in \S2. In either case, our key observation is the
following transformation law of $\log t$, where $t$ is a conformal
scale (a fiber coordinate of the bundle $\G$ determined by a section
$g\in[g]$):
$$
  -\log\wh t=-\log t+\Up \quad\text{whenever}\quad\wh g=e^{2\Up}g.
$$
In view of (a), we extend $t$ off $\G$ and apply $\wt\Delta^{n/2}$ to
$-\log t$. Then we see that its restriction to $\G$ gives the $Q$-curvature;
the required transformation law is clear from that of $-\log t$. This
formulation is naturally related to Branson's one. In the argument of
analytic continuation in dimension $m$ the function $\log t$ appears as
the differential in $m$ of the density $t^{(n-m)/2}$ at $m=n$.
Corresponding to (b), we also express the $Q$-curvature as an obstruction
to the existence of a smooth solution to $\wt\Delta F=0$ with an initial
condition $F|_\G=\log t$.  Since the Laplace equation in the ambient metric
can be reformulated as the one in the Poincar\'e metric, we see
that this derivation is equivalent to that of Fefferman-Graham, mentioned above.

Note that there is another derivation of the invariant operators in terms
of a bundle calculus associated with the conformal Cartan connection, which
is called tractor calculus \cite{CG}. Corresponding to this derivation,
Gover and Peterson \cite{GP} gave a tractor expression of the $Q$-curvature.
Gover informed us that their formula can be translated into an ambient metric
expression that is equivalent to our construction corresponding to (a).
See also the remark at the end of \S2.

In \S3, we turn to CR geometry. For strictly pseudoconvex CR structures,
the ambient metric is given as a Lorentz K\"ahler metric. The powers of the
ambient Laplacian induce invariant powers of the sublaplacian $P_{n/2}$ and the
ambient construction of $Q$-curvature is also valid. (For a comprehensive
treatment of CR invariant operators see \cite{GG}.)  A new feature in this
setting is that $P_{n/2}$ has a large null space, including the space of CR
pluriharmonic functions $\mathcal P$.  Thus the $Q$-curvature on
$(2N-1)$-dimensional CR manifolds, $Q_{\theta}^\CR$, which is a local invariant
of the pseudohermitian structure $\theta$, satisfies
\begin{equation}\label{CR-trans}
 e^{2N\up} Q^\CR_{\wh\theta}=Q^\CR_{\theta}\quad
 \text{whenever $\wh\theta=e^{2\up}\theta$ with $\up\in\mathcal P$}.
\end{equation}
If $N=2$, it has been shown in \cite{H} that this transformation law uniquely
characterizes $Q^\CR_{\theta}$ up to a constant multiple. As a consequence, we see
that the leading term of the logarithmic singularity of the Szeg\"o kernel,
$\psi$, is a constant multiple of $Q^\CR$, since the Szeg\"o kernel enjoys
the same transformation law.  While such a simple characterization does not
hold for higher dimensions, the transformation law still indicates an intimate
link between the $Q$-curvature and the Szeg\"o kernel.

\section{$Q$-curvature in conformal geometry}

\subsection{Conformally invariant operators}
We first recall basic materials on the ambient metric construction of the
invariant operators from \cite{FG1} and \cite{GJMS}.

Let $(M,[g])$ be a conformal manifold of signature $(p,q)$, $p+q=n\ge3$.
Then $M$ admits the metric bundle ${\mathcal G}\subset S^2 T^*M$, a ray
bundle consisting of the metrics in the conformal class $[g]$. There are
dilations $\delta_s:\G\to\G$ given by $\delta_s(g)=s^2g$ for $s>0$, and
the homogeneous functions, with respect to $\delta_s$, on $\G$ are
called conformal densities; the space of densities of weight $w$ is denoted
by $\E(w)$, i.e.,
$$
\E(w)=\{f\in C^\infty(\G): \delta_s^*f=s^wf\text{ for any } s>0\}.
$$
Conformally invariant operators are then defined as operators acting on
the conformal densities:
$$
P:\E(w)\to\E(w').
$$
A choice of representative $g\in[g]$ determines a trivialization
\begin{equation*}\label{E-trivialization}
 \E(w)\ni f\mapsto f_g:=f\circ g\in C^\infty(M)
\end{equation*}
such that $f_{\wh g}=e^{w\Up}f_g$ when $\wh g=e^{2\up}g$. Thus an invariant
operator $P$ defines, for each representative $g\in[g]$, an operator $P_g$
on $C^\infty(M)$ such that
\begin{equation*}\label{trans-P}
 P_{\wh g}=e^{w'\up}P_g e^{-w\up}
 \quad\text{whenever } \wh g=e^{2\up}g.
\end{equation*}
In particular, if $w=0$, then $\E(0)=C^\infty(M)$ and $P$ acts on the
functions on $M$.  We say that $P$ is an invariant diferential operator if
each $P_g$ is given by a differential operator on $C^\infty(M)$.

The ambient metric $\wt g$ is formally defined on  $\wt\G=\G\times(-1,1)$
along $\G$ which is now embedded as a hypersurface
$\iota:\G\to\G\times\{0\}\subset\wt\G$. It is characterized by the
following three conditions:

\smallskip

(1) $\wt g$ is an extension of the tautological two-tensor $g_0$ on $\G$,
i.e.,  $\iota^*\wt g=g_0$;

(2) $\delta_s^*\wt g=s^2 \wt g$ for any $s>0$;

(3) $\wt g$ is an asymptotic solution to $\operatorname{Ric}(\wt g)=0$
along $\G$.

\smallskip

\noindent
When $n$ is odd, these conditions uniquely determine a formal power series
of $\wt g$ up to homogeneous diffeomorphisms that fix $\G$; but when $n$ is
even, $\wt g$ exists in general only to order $n/2$.

Many invariant differential operators can be constructed out of the
ambient metric. The basic procedure is to construct a differential
operator on the ambient space $\wt\G$ which preserves the homogeneity
$\wt P:\wt\E(w)\to\wt\E(w')$ and then prove that $\wt P$ induces an
operators $P:\E(w)\to\E(w')$, namely, prove that $(\wt P\wt f)|_\G$
depends only on $\wt f|_\G$. We then call $\wt P$ an ambient extension
of $P$.  For example, the powers of the Laplacian
$\wt\Delta=-\wt\nabla_I\wt\nabla^I$ in the ambient metric
$$
\wt P_k=\wt\Delta^{k}:\wt\E(k-n/2)\to\wt\E(-k-n/2)
$$
for $k>0$ (and $k\le n/2$ if $n$ is even) induce
$$
P_k:\E(k-n/2)\to\E(-k-n/2).
$$
The leading part of $P_k$ for each representative $g$ is the powers of the
Laplacian $\Delta^k$ in $g$ and hence $P_k$ is called an invariant power
of the Laplacian -- see Proposition 2.1 of \cite{GJMS}.

More operators have been constructed by Alexakis [1] by using the
harmonic extension of densities. Denoting by $T$ the infinitesimal
generator of the dilations $T=\frac{d}{ds}\delta_s|_{s=1}$,
we set $\wt\rho=\|T\|^2$; then $\wt\rho\in\E(2)$ and $\wt\rho=0$
defines $\G$. Then the harmonic extension of densities are
explicitly given by the following lemma, which is
a part of Proposition 2.2 of \cite{GJMS}.

\begin{lemma}
Let $f\in\E(w)$ and set $k=n/2+w$. If $k\not\in\{1,2,\dots\}$, then
$f$ admits an extension to $\wt f_m\in\wt\E(w)$ such that
$\wt\Delta\wt f_m=O(\wt\rho^m)$ for any $m\ge0$. Such an $\wt f_m$ is
unique modulo $O(\wt\rho^{m+1})$ and is given by
$$
 \wt f_m=E_mE_{m-1}\cdots E_1\wt f,
 \quad\text{where}\quad
 E_l=1+\frac{1}{4l(k-l)}\wt\rho\,\wt\Delta
$$
and $\wt f\in\wt\E(w)$ is an arbitrary extension of $f$. If
$k\in\{1,2,\dots\}$, then the same result is true
with the restriction $m<k$.
\end{lemma}

Using this harmonic extension and $\wt\nabla^p\wt R$, the iterated
covariant derivative of the curvature tensor of the ambient metric
$\wt g$, we form a complete contraction
$$
 \wt P\wt f=\contr\Big(\wt\nabla^{p_1}\wt R\otimes\cdots\otimes
 \wt\nabla^{p_l}
 \wt R\otimes\wt\nabla^q\wt f_m\Big).
$$
It defines a map $\wt P:\wt\E(w)\to\wt\E(w')$, where
$w'=w-p_1-\cdots-p_l-2l-q$. If $q$ is sufficiently small (e.g., $q\le m$)
then the lemma above ensures that $\wt P$ induces an invariant operator
$P:\E(w)\to\E(w')$.

\medskip
{\em Remark.}
In \cite{A}  it is shown that all conformally invariant differential
operators $P:\E(w)\to\E(w')$ arise as above provided the dimension $n$
is odd and $n/2+w\not\in\{1,2,\dots\}$. The result of \cite{A} applies
also to nonlinear operators.

\subsection{$Q$-curvatures in terms of the ambient metric}

We now define $Q$-curvatures for the invariant operators constructed
as above.

\begin{theorem}\label{thm-Q-a}
Let $P:\E(0)\to\E(w)$ be an invariant differential operator with an
ambient extension $\wt P:\wt\E(0)\to\wt\E(w)$.  For $g\in[g]$, choose
$t\in\wt\E(1)$ such that $t(u^2g)=u$ on
$\G$. Then
$$
 Q_g=-\big(\wt P \log t\big)\circ g
$$
is independent of the extension of\/ $t$ off $\G$ and defines a function
determined by $g$. Moreover, if $P1=0$, then $Q_g$ satisfies the
transformation law
\begin{equation}\label{trans-Q}
 e^{-w\Up}Q_{\wh g}=Q_g+P_g\Up\quad
 \text{ whenever }\quad \wh g=e^{2\Up} g.
\end{equation}
\end{theorem}

{\em Proof.}\/
If $t'$ is another function $\wt\E(1)$ which agrees with
$t$ on $\G$, then we have
$t'=e^{f}t$ for an $f\in\wt\E(0)$ such that $f=O(\wt\rho)$. So
$$
\wt P\log t'-\wt P\log t=\wt Pf=O(\wt\rho)
$$
and hence $(\wt P\log t')\circ g=(\wt P\log t)\circ g$. To prove
the transformation law, we  extend  $\Up$ to $\Up\in\wt\E(0)$ and
set $\wh t=e^{-\Up} t$. Since $\wt P$ has no zeroth order term, we have
$\wt P\log t\in\wt\E(w)$ so that
$$
e^{-w\Up}Q_{\wh g}=-e^{-w\Up}\big(\wt P \log \wh t\big)\circ{\wh g}
=-\big(\wt P \log\wh t)\circ g.
$$
Substituting $\wt P\log \wh t=\wt P\log t-\wt P\Up$ into the right-hand
side, we get \eqref{trans-Q}.
\qed

\medskip

In particular, if $P$ is the invariant power of the Laplacian $P_{n/2}$
then $w=-n$ and $Q_g=-\big(\wt\Delta^{n/2}\log t\big)\circ g$. We now
show that this $Q_g$ agrees with the $Q$-curvature defined by Branson
\cite{Br}, which we recall briefly. For $m\ge n/2$, we denote by
$P_{n/2,m}$ the invariant powers of Laplacian of order $n$ in dimension
$m$. Let $\wt Q_{n/2,m}$ be the zeroth order term of $P_{n/2,m}$ in
the metric $g$. Then, noting $P_{n/2,n}1=0$, we may write
$\wt Q_{n/2,m}=(m-n)/2\, Q_{n/2,m}$ for a scalar Riemannian invariant
of $g$. Moreover, {from} the construction of $P_{n/2,m}$,
we see that $Q_{n/2,m}$ is expressed as a linear combination complete
contractions of the tensor products of $\nabla^l R$, the coefficients
of which are rational in $m$ and regular at $m=n$. Thus we may
substitute $m=n/2$ and define $Q$-curvature by $Q_{n/2,n}$.

In the identification \eqref{E-trivialization} with respect to $g$, the
constant function $1\in C^\infty(M)$ corresponds to $t^w\in\E(w)$.
Thus, extending $t^w$ to $\wt\E(w)$, we have
$$
\wt Q_{n/2,m}=P_{n/2,m}1=\big(\wt\Delta^{n/2}t^{(n-m)/2}\big)\circ g.
$$
Substituting $t^{(n-m)/2}=1+(n-m)/2\,\log t+O\big((n-m)^2\big)$
into the right-hand side gives
$$
 \wt Q_{n/2,m}=\frac{(n-m)}{2}\big(\wt\Delta^{n/2}\log t\big)\circ g
 +O\big((n-m)^2\big),
$$
which implies $Q_{n/2,n}=-\big(\wt\Delta^{n/2}\log t\big)\circ g$
as claimed.

\subsection{$Q$-curvature in terms of Poincar\'e metrics}
There is another derivation of $P_{n/2}$, which was also given in
\cite{GJMS}. For $f\in\E(0)$, take an extension $\wt f\in\wt\E(0)$
such that $\wt\Delta\wt f=O(\wt\rho^{n/2})$ -- see Lemma 2.1. Then
$\wt\rho^{\,1-n/2}\wt\Delta\wt f|_\G$ is shown to agree with
$c_n P_{n/2}f$, where $c_n=2^{2-n}((n/2-1)!)^{-2}$.
This derivation of $P_{n/2}$ can be reformulated as follows:

\begin{lemma} \label{Laplace-log}
Let $n$ be even. For a representative $g\in[g]$, take $t$ as in Theorem
$\ref{thm-Q-a}$ and set $\rho=\wt\rho/(2t^2)$. Then, for each $f\in\E(0)$,
there exists a formal solution to $\wt\Delta\wt f=0$ of the form
$\wt f=\wt f_0+\eta \,\wt\rho^{\,n/2}\,\log \rho$ with
$\wt f_0\in\wt\E(0)$ such that $\wt f_0|_\G=f$ and $\eta\in\wt\E(-n)$.
Here $\wt f_0\mod O(\rho^{n/2})$ and $\eta\mod O(\rho^\infty)$ are
determined by $f$, and moreover, $\eta|_\G$ is a non-zero constant
multiple of $P_{n/2}f$.
\end{lemma}

{\em Proof.}
By taking the smooth part and the log term, we decompose
$\wt\Delta\wt f=0$ into a system of equations
$$
\wt\Delta\wt f_0+[\wt\Delta,\log \rho]\,\eta\wt\rho^{\,n/2}=0,\quad
 \wt\Delta\,\eta=0.
$$
We solve this system by using Lemma 2.1. Noting that
$$
[\wt\Delta,\log \rho]\,\eta\wt\rho^{\,n/2}
=2n\,\eta\wt\rho^{\,n/2-1}+O(\rho^{n/2}),
$$
we first solve $\wt\Delta\wt f_0=O(\rho^{\,n/2-1})$ and $\wt\Delta\eta=0$
under the initial conditions $\wt f_0|_\G=f$ and
$\eta|_\G=-c_n/(2n)\, P_{n/2}f$. Then we have
$\wt\Delta\wt f_0+[\wt\Delta,\log \rho]\,\eta\wt\rho^{\,n/2}
=O(\rho^{\,n/2})$, and thus we may modify $\wt f_0$ so that
$\wt\Delta\wt f_0+[\wt\Delta,\log \rho]\,\eta\wt\rho^{\,n/2}=0$.
The uniqueness of $\eta$ is clear from this construction. \qed

\medskip

Corresponding to this derivation of $P_{n/2}$, we have the following
characterization of the $Q$-curvature.

\begin{theorem}\label{thm-log-solution}
Let $n$ be even. For a representative $g\in[g]$, take $t$ as in
Theorem $\ref{thm-Q-a}$ and set $\rho=\wt\rho/(2t^2)$. Then there
is a formal solution to $\wt\Delta F=0$ of the form
$$
 F=\log t+\varphi+\eta\, \wt\rho^{n/2}\log\rho
$$
with
$\varphi\in \wt\E(0)$ such that $\varphi=O(\rho)$ and $\eta\in \wt\E(-n)$.
Here $\varphi\mod O(\rho^{n/2})$ and $\eta\mod O(\rho^\infty)$ are
determined by $g$ and, moreover, $\eta\circ g$ is a constant multiple of
the $Q$-curvature of $P_{n/2}$.
\end{theorem}

The proof of this theorem is just a straightforward modification of
that of Lemma \ref{Laplace-log}; the last statement follows form the
fact that $\eta$ is a multiple of $\wt\Delta^{n/2}\log t$ on $\G$.
We will omit the details and, instead, we show that this theorem is
equivalent to Theorem 3.1 of \cite{FG2}, which we state as
Theorem \ref{Thm-FG} below.

Let $X=M\times(0,1)$ and identify $M$ with a portion of the boundary
$M\times\{0\}$.  The Poincar\'e metric $g_+$ is a metric on $X$
satisfying the following conditions: $g_+$ satisfies the Einstein equation
$\operatorname{Ric}(g_+)+ng_+=0$ asymptotically along $M$, and if $r$ is a
defining function of $M$ in $X$, then $h=r^{-2}g_+$ is smooth on
$\overline X=M\times [0,1]$ and $h|_{TM}\in [g]$. Note that $r\mod O(r^2)$
corresponds to a representative $g\in[g]$. The higher jets of $r$ can be
uniquely determined by the normalization $\|d\log r\|_{g_+}=1$.

\medskip
\noindent
\begin{theorem} {\rm (\cite{FG2})}
\label{Thm-FG}
Let $n$ be even. For a representative $g\in[g]$, take a defining
function $r$ such that $\|d\log r\|_{g_+}=1$ and $(r^2g_+)|_{TM}=g$.
Then, there is an asymptotic solution to the equation
$$
 \Delta_+ U=n
$$
of the form
$$
 U=\log r+A+B\,r^n\log r
$$
with $A, B\in C^\infty(\overline X)$ which are even in $r$ and $A|_M=0$. Here
$A\mod O(r^n)$ and $B\mod O(r^\infty)$ are formally determined by $g$,
and moreover, $B|_M$ is a constant multiple of the $Q$-curvature.
\end{theorem}

To translate this theorem into Theorem \ref{thm-log-solution}, we recall
the relation between the ambient metric and Poincar\'e metric from
\cite{FG1} and \cite{GZ}.  With respect to a suitable a decomposition
of $\wt\G=\R_+\times M\times (-1,1)$, we have
\begin{equation}\label{normal-g}
  \wt g=2tdtd\rho+2\rho dt^2+t^2 \overline g_\rho,
\end{equation}
where $t\in\R_+$ is homogeneous of degree $1$, $\rho\in(-1,1)$ and
$\overline g_\rho$ is a one-parameter family of metrics on $M$ such that
$\overline g_0=g$. Then $T=t\partial_t$ and $\wt\rho=2t^2\rho$ hold.
We embed $\overline X=M\times [0,1]\ni(x,r)$ into $\wt\G$ by the map
$\iota(x,r)=(1/r,x,-r^2/2)$ so that $\iota(\overline X)=\{2\wt\rho=-1,\rho\le0\}$.
Then $g_+=\iota^*\wt g$ gives the Poincar\'e metric. Now set
$s=t\sqrt{-2\rho}=t\,r$ and define new coordinates $(s,x,r)$ of $\wt\G$ in
 which $\iota(X)=\{s=1\}$. Then we have $\wt g=s^2g_+-ds^2$ and hence
$$
\wt\Delta=s^{-2}(\Delta_++(s\pa_s)^2+ns\pa_s),
$$
where $\Delta_+$ is considered as an operator in the variables $(x,r)$.
Thus
$$
\begin{aligned}
 s^2\wt\Delta F&=
 s^2\wt\Delta(\log s-\log r+\varphi+\eta\,\wt\rho^{\,n/2}\log\rho)\\
& =n-\Delta_+(\log r-\varphi-\eta\,\wt\rho^{\,n/2}\log\rho)
\end{aligned}
$$
because $-\log r+\varphi+\eta\,\wt\rho^{\,n/2}\log\rho$ is homogeneous of
degree $0$. Now, restricting the both sides to $s=1$, we get
$$
 (\wt\Delta F\big)\big|_{s=1}=n-\Delta_+ U,
$$
where  $U=-F\big|_{s=1}$ and it is of the form
$$
 U=\log r+A+B\,r^n \log r.
$$
Therefore $\wt\Delta F=0$ is equivalent to $\Delta_+ U=n$. Comparing
the log term coefficients, we have
$B(x,r)=-2(-1)^{n/2}\eta|_{(t,x,\rho)=(1,x,-r^2/2)}$.

\subsection{Examples}
We  give two examples of the pairs $(P,Q)$, an invariant operator $P$ and
the associated $Q$-curvature. It is a routine computation and we only outline
the computation by quoting basic formulas from \cite{FG1} and \cite{GJMS}.

Fixing a representative $g\in[g]$, we take local coordinates $(t,x_i,\rho)$
of $\wt\G=\R_+\times M\times (-1,1)$ such that \eqref{normal-g} holds;
we here rename the coordinates as $(x_I)=(x_0,x_i,x_\infty)$ and use capital
indices $I,J,K,\dots$ (resp. small indices $i,j,k,\dots$) to ran through
$0,1,\dots,n,\infty$ (resp. $1,\dots, n$). With these coordinates, it is
easy to compute the covariant derivatives of $\log t$. We have
$\wt\nabla_I\wt\nabla_J\log t=0$ except for the following two cases:
\begin{equation}\label{nabla-log}
 \wt\nabla_0\wt\nabla_0\log t=-t^{-2},
 \quad\wt\nabla_i\wt\nabla_j\log t=t^{-2}P_{ij}+O(\rho).
\end{equation}
Here $P_{ij}$ is the Rho tensor, the trace modification of the Ricci tensor
of $g$, determined by
$$
 P_{ij}=\frac{1}{n-2}(R_{kij}{}^k-J\,g_{ij}),
 \quad J=P_i{}^i=\frac{1}{2(n-1)}R_{ij}{}^{ji}.
$$
In particular, we see that $\wt\Delta \log t|_{\G}=-J$ is a constant
multiple of the scalar curvature.

We next express the components of $\wt R_{IJKL}$ on $\G$ in terms of
the curvature of $g$:
\begin{equation}\label{wtR}
 \wt R_{IJK0}=0,\
 \wt R_{ijkl}=t^2W_{ijkl},\
 \wt R_{ijk\infty}=t^2C_{kij},\
 \wt R_{\infty ij \infty}=\frac{t^2 B_{ij}}{n-4}.
\end{equation}
(If $n=4$, $\wt R_{\infty ij \infty}$ is undetermined.) Here $W$, $C$
and $B$, called the Weyl, Cotton and Bach tensor, respectively, are defined
as follows: $W_{ijkl}$ is the totally trace-free part of the curvature
tensor $R_{ijkl}$; $C_{ijk}=\nabla_k P_{ij}-\nabla_j P_{ik}$,
and $B_{ij}=\nabla^k C_{ijk}+P^{kl}W_{kijl}$. These relations and the usual
symmetries of the curvature tensor determine all the components of $\wt R$.

Our first example of $P$ is the operator $P^1:\E(0)\to\E(-6)$ with the ambient
expression
$$
 \wt P^1(\wt f)=\wt R_{IJK}{}^L\wt R^{IJKM}\wt\nabla_L\wt\nabla_M \wt f_1.
$$
Recalling $\wt f_1=\wt f+\frac{1}{2(n-2)}\wt\rho\,\wt\Delta\wt f$ and using
\begin{equation}\label{GJMS-eq}
 \wt R_{IJK}{}^LT_L=0,\quad
 \wt\nabla_IT_J=g_{IJ}
 \text{ with} \quad 2\,T_I=\wt\nabla_{\!I}\,\wt\rho
\end{equation}
 (see  (1.3)--(1.6) of \cite{GJMS}), we see that
$$
\wt P^1(\wt f)=\wt R_{IJK}{}^L\wt R^{IJKM}\wt\nabla_L\wt\nabla_M \wt f
 +\frac{1}{n-2}\|\wt R\|^2\wt\Delta \wt f+O(\rho).
$$
In terms of the representative metric $g$, it can be expressed as
$$
 P^1(f)=
 W_{ijk}{}^l W^{ijkm}\nabla_l\nabla_mf
 -2C_{kij}W^{ijkl}\nabla_lf+\frac{1}{n-2}\|W\|^2\Delta f.
$$
On the other hand, using \eqref{nabla-log} and \eqref{wtR}, we can express
the $Q$-curvature $ Q^1=-\wt P^1(\log t)\circ g$ as
$$
 Q^1 =-W_{ijk}{}^l W^{ijkm}P_{lm}+\|C\|^2+\frac{1}{n-2}\|W\|^2J.
$$

Our next example is $P^2:\E(0)\to\E(-6)$ induced by
$$
 \wt P^2(\wt f)=\wt R_{IJKL}(\wt\nabla^M\wt R^{IJKL})\wt\nabla_M \wt f_1.
$$
Noting that $2\wt P^2(\wt f)=\wt\nabla^I(\|\wt R\|^2\wt\nabla_I \wt f_1)$,
we have
$$
 2\wt P^2(\wt f)=\wt\nabla^I\big(\|\wt R\|^2\wt\nabla_I \wt f\big)
 +\frac{n-6}{n-2}\|\wt R\|^2\wt\Delta \wt f+O(\rho)
$$
and from this one can easily deduce
$$
 2P^2(f)=\nabla^i\big(\|W\|^2\nabla_if)+\frac{n-6}{n-2}\|W\|^2\Delta f.
$$
The $Q$-curvature $Q^2=-\wt P^2(\log t)\circ g$ is
$$
 Q^2
 =2\Big(W_{ijkl}\nabla^iC^{jkl}+W_{ijk}{}^lW^{ijkm}P_{lm}
 +2\|C\|^2-\frac{1}{n-2}\|W\|^2J\Big).
$$

\medskip

{\em Remark.}
The $Q$-curvature associated with $P^2:\E(0)\to\E(-6)$ in dimension $6$
has been obtained in Gover-Peterson \cite{GP}; their $Q$-curvature is
$\frac{1}{8}\Delta\|W\|^2$. Our $Q^2$ is consistent with their formula because
$Q^2=\frac{1}{8}(\Delta\|W\|^2-\wt\Delta\|\wt R\|^2\circ g)$ when $n=6$
and $\wt\Delta\|\wt R\|^2\circ g$ is a conformal invariant -- see \cite{FG1}.

\section{$Q$-curvature in CR geometry}
\subsection{Ambient metric and invariant contact forms}
We now turn to CR geometry.  We first recall the ambient metric of \cite{F1}
and \cite{F2}. Let $M$ be a strictly pseudoconvex real hypersurface in $\C^N$
and let $J$ be the complex \MA operator
\begin{equation}\label{MAop}
J[\rho]=(-1)^N \det\left(
      \begin{array}{cc}
        \rho&  \rho_j\\
        \rho_{\overline k}&  \rho_{j\overline k}
      \end{array}
    \right),
\quad  \rho_{j}=\frac{\pa\rho}{\pa{z_j}}, \text{ etc.}
\end{equation}
Then there is a smooth defining function of $M$ that is positive on the
pseudoconvex side and satisfies $J[u]=1+O(u^{N+1})$; such a $u$ is unique
modulo $O(u^{N+2})$. The ambient metric lives on $\C^*\times \wt M$
for a small collar neighborhood $\wt M$ of $M$. It is the Lorentz-K\"ahler
metric
\begin{equation}\label{gCR}
  \wt g[u]=-\sum_{j,k=0}^{N}
\frac{\pa^2\big(|z_0|^2u(z)\big)}{\pa z_j\pa \overline z_k}\,dz_j d\overline z_k,
\end{equation}
with $(z_0,z)\in\C^*\times \wt M$. Note that $J[u]=1+O(u^{N+1})$ implies
Ric$(\wt g)=O(u^N)$.

The defining function $u$ also specifies a contact form
$\theta[u]=\Im\, \pa u|_{TM}$ of $M$ and $\wt g[u]$ induces a real
Lorentz metric $g[u]$ on the circle bundle $S^1\times M$. Since $g[u]$ is
shown to depend only on $\theta[u]$, we may write the metric as $g[\theta]$.
This correspondence $\theta\mapsto g[\theta]$ can be extended to a general contact
form $\theta$ of $M$ in such a way that $g[e^{2\Up} \theta]=e^{2\Up}g[\theta]$ holds,
and we have a conformal class of Lorentz metric $[g]$ on $S^1\times M$ --
see \cite{L0}. For the conformal manifold $(S^1\times M,[g])$, the metric
bundle and the ambient space are given by $\G=\C^*\times M$ and
$\wt\G=\C^*\times\wt M$ respectively, and the metric $\wt g[u]$ satisfies
the conditions (1), (2) and (3) of \S2.1. Thus the definition of the ambient
metric in conformal and CR cases are compatible, where $(2N-1)$-dimensional
CR manifolds correspond to $2N$-dimensional Lorentzian conformal manifolds
-- see \cite{FG1}.

The contact form $\theta[u]$ defined above has special importance, and we call
$\theta[u]$ an invariant contact form.  This notion can be generalized to CR
manifolds $M$: $\theta$ is an invariant contact form on $M$ if it is locally
given as an invariant contact form for some local embedding of $M$ into $\C^N$.
An intrinsic formulation of invariant contact form is also given (\cite{Fa},
\cite{L}): $\theta$ is an invariant contact form if it is locally
volume-normalized with respect to a closed $(N,0)$-form on $M$. {From} this
characterization, it is straightforward to see that any two invariant
contact forms $\theta$ and $\theta'$ satisfy $\theta'=e^{2\up}\theta$ with a CR
pluriharmonic function $\up$ (that is, $\up$ is locally the real part of a
CR function). Note that, when $N\ge3$, Lee \cite{L} showed that $\theta$ is an
invariant contact form if and only if $\theta$ is pseudo-Einstein, that is,
the Tanaka-Webster Ricci tensor of $\theta$ is a scalar multiple of the Levi
form (this condition is vacuous when $N=2$).

We next consider the CR analog of the operators $P_k$. CR densities of weight
$(w,w)$ are functions $f(z_0,z)$ on $\G$ such that
$f(\lambda z_0,z)=|\lambda|^{2w}f(z_0,z)$ for any $\lambda\in\C^*$. The
totality of such functions is denoted by $\E(w,w)$. For each $\theta$, the metric
$g[\theta]$ determines a $S^1$-subbundle of $\pi:\G\to M$. Restricting each
$f\in\E(w,w)$ to the circle bundle, we obtain a function $\pi_* f$ on $M$;
this correspondence gives an identification $\E(w,w)\cong C^\infty (M)$.
Note also that $\E(w,w)$ can be regarded as a subspace of conformal densities
$\E(2w)$ for the conformal manifold $(S^1\times M,[g])$.  As in the conformal
case, we extend $\E(w,w)$ to the ambient space and define $\wt\E(w,w)$ to be
the smooth functions on $\wt\G$ which are homogeneous of degree $(w,w)$ in
$z_0$ variable. Then the powers of the ambient Laplacian $\wt\Delta^k$ maps
$\wt\E(w,w)$ into $\wt\E(w-k,w-k)$ and, for $w=k-N\le0$, it induces an operator
$P_k:\E(w,w)\to\E(w-k,w-k)$.
{From} this construction, it is clear that the CR invariant operator $P_k$
is the restriction of the conformally invariant operator
$P_k:\E(2w)\to\E(2w-2k)$.

\subsection{CR $Q$-curvature}
 The CR version of the $Q$-curvature is defined by
$$
Q^\CR_\theta:=\pi_*Q_g,
$$
where $Q_g$ is the $Q$-curvature, of the conformal $P_N$, in the metric
$g=g[\theta]$, and where $\pi$ is the projection $S^1\times M\to M$. Since $Q_g$
is $S^1$-invariant and pushes forward to a function on $M$. Then, as in
the conformal case, we have
\begin{equation}\label{CRQ-trans}
 e^{2N\up}Q^\CR_{\wh\theta}=Q^\CR_\theta+P_N\up
 \quad\text{whenever }\wh\theta=e^{2\up}\theta.
\end{equation}
Here $P_N$ is computed in $\theta$.

\begin{proposition}\label{prop-Q-vanishing}
If $\theta$ is an invariant contact form, then $Q^\CR_\theta=0$.
\end{proposition}

{\em Proof.}
For the metric $g[\theta]$, we may take $|z_0|^2$ as a fiber coordinate of
$\C^*\times M\to S^1\times M$. Then
$Q^\CR_\theta=-(\wt\Delta^N\log |z_0|^2)\circ g[\theta]=0$ because $\wt\Delta$
kills pluriharmonic functions.
\qed

\medskip

In view of this proposition, we have another expression of $Q^\CR$.
Take an invariant contact form $\theta_0$ as a reference and set
\begin{equation}\label{CR-Q-P}
  Q^\CR_\theta=e^{-2N\up}P_N\up,
\end{equation}
where $\theta=e^{2\up}\theta_0$ and $P_N$ is computed in $\theta_0$.
This is well-defined because $\up$ modulo additions of CR pluriharmonic
functions is independent of the choice of $\theta_0$ and CR pluriharmonic
functions are killed by $P_N$.

If $M$ is a real hypersurface in $\C^N$,
then $M$ admits a global invariant contact form $\theta$ so that
$Q^\CR_\theta=0$ on $M$. However, for abstract CR manifolds $M$, there is
a topological obstruction for the global existence of an invariant
contact form $\theta$: the existence of $\theta$ implies the vanishing of
the first Chern class of the holomorphic tangent bundle
$c_1(T^{1,0}M)$ in $H^2(M,\R)$.
This obstruction, for $N\ge3$, was first found by Lee \cite{L} in the study
of pseudo-Einstein contact form and his argument implicitly contains the
proof for the $N=2$ case. At present, we do not know if we can always
choose $\theta$ so that $Q^\CR_\theta$ vanishes globally.

\medskip

\noindent
{\em Remark.}
The operators $P_N$ were first introduced in Graham \cite{G0} as a
compatibility operator for the Dirichlet problem for the Bergman
Laplacian for the ball in $\C^N$. This construction of $P_N$ was
generalized, in \cite{GL}, to the boundaries of strictly pseudoconvex
domains in $\C^N$; the CR invariance of $P_2$ was realized later in \cite{H}.
Graham \cite{G0} proved that the kernel of $P_N$ for the sphere agrees
with the space of CR pluriharmonic functions (this can be partially
generalized to the curved case \cite{GL}). As a result, on the sphere,
we see from the expression \eqref{CR-Q-P} that $Q_\theta^\CR=0$ if and only
if $\theta$ is an invariant contact form.

\subsection{Logarithmic singularity of the Szeg\"o kernel}
Now let $M$ be the boundary of a strictly pseudoconvex domain
$\Omega$ in $\C^N$. For a choice of contact form $\theta$ on $M$,
we define $H_\theta^2(M)$ to be the kernel of $\overline\pa_b$
in $L^2(M)$ with respect to the volume element $\theta\wedge (d\theta)^{N-1}$.
Then the Szeg\"o kernel $K(x,y)$ is defined as the reproducing kernel of
the Hilbert space $H^2_\theta(M)$. $K(x,y)$ can be extended to a holomorphic
function $K(z,\overline w)$ on $\Omega\times\overline\Omega$; its restriction to the
diagonal $K(z,\overline z)$ admits an expansion
$$
 K(z,z)=\varphi(z)u(z)^{-N}+\psi(z)\log u(z),
$$
where $\varphi$, $\psi$ are functions smooth up to the boundary and $u$
is a defining function of $\Omega$ -- see \cite{F2}, \cite{BS}. This asymptotic
expansion is locally determined by the CR structure of $\pa\Omega$ and $\theta$.
Moreover, $\psi_\theta=\psi|_M$ is shown to be a local pseudohermitian invariant
of $\theta$, that is, $\psi_\theta$ can be written as a linear combination of
complete contractions, with respect to the Levi form, of the tensor products
of Tanaka-Webster curvature and torsion and their covariant derivatives.

In general, there is no simple transformation law of the Szeg\"o kernel
under the scaling of contact form $\wh\theta=e^{2\up}\theta$. But if $\up$ is CR
pluriharmonic, we have $\wh K=e^{-2N\up}K$, where $\up$ is extended to a
pluriharmonic function in $\Omega$. In particular, we obtain
$$
 e^{2N\up}\psi_{\wh\theta}=\psi_\theta.
$$
In case $N=2$, this transformation law is strong enough to characterize
$\psi_\theta$ up to a constant multiple. In fact, we have

\medskip
\noindent{\bf Theorem.}{ (\cite{H}) }{\em
Let $N=2$. Suppose that $S_\theta$ is a scalar  pseudohermitian invariant
satisfying the transformation law
\begin{equation}\label{transCRP}
  e^{2N\up}S_{\wh\theta}=S_\theta
  \quad\text{whenever $\wh\theta=e^{2\up}\theta$
   and $\up$ is CR pluriharmonic.}
\end{equation}
Then $S_\theta$ is a constant multiple of
\begin{equation}\label{Q-formula}
\Delta_b R-2\,\Im
\nabla^\alpha\nabla^\beta A_{\alpha\beta},
\end{equation}
where $R$ is the Tanaka-Webster scalar curvature, $A$ is the torsion,
$\Delta_b$ is the sublaplacian computed in $\theta$.
}
\medskip

Since $Q^\CR_\theta$ also satisfy the transformation law \eqref{transCRP},
the theorem above implies $\psi_\theta=c\, Q_\theta^\CR$ for a universal
constant $c$, which can be identified by an explicit computation for
an example (c.f. \cite{H}, \cite{GG}). Thus we have

\begin{proposition}\label{Szego}
If $N=2$, then
$$
32\pi^2\psi_\theta=
Q_\theta^\CR=\frac{4}{3}\big(\Delta_b R-2\,\Im
\nabla^\alpha\nabla^\beta A_{\alpha\beta}\big).
$$
\end{proposition}

For $N\ge3$, there are examples of pseudohermitian invariants that
satisfy \eqref{transCRP} for any $\up\in C^\infty(M)$ -- see \cite{F2}.
Such invariants are called CR invariants of weight $N$. Thus it is a
natural conjecture that $S$ is a constant multiple of $Q^\CR$ up to an
addition of CR invariant of weight $N$. (In case $N=2$, there is no
CR invariant of weight $2$ -- see \cite{Gr1}, and this conjecture is reduced
to the theorem above.)

We finally note that the integral of the $Q$-curvature
$L_M=\int_M Q_\theta^\CR\theta\wedge(d\theta)^{N-1}$ is independent of the choice
of a contact form $\theta$ and gives a CR invariant; this follows from the
analogous fact in the conformal case.  In case $N=2$, it turns out that
$L_M=0$ because \eqref{Q-formula}  is the divergence of the one form
$\nabla_\alpha R-i\nabla^\beta A_{\alpha\beta}$. We also see from the
argument of \S3.2 above that $L_M$ vanishes if $M$ admits a global
invariant (or pseudo-Einstein) contact form. It should be interesting to
find a link between $L_M$ and the Chern class $c_1(T^{1,0}M)$,
which obstructs the existence of an invariant contact form $\theta$.

\end{document}